\documentclass[12pt]{article}
\usepackage{geometry}
\usepackage{chngpage}
\usepackage{graphicx}
\usepackage{amsmath}
\usepackage{amsfonts}
\usepackage{amssymb}
\usepackage{latexsym}
\usepackage[brazil]{varioref}
\usepackage[english,brazil]{babel}
\usepackage[latin1]{inputenc}
\usepackage[dvips]{color}

\newtheorem{defn}{Defini\c{c}{\~a}o}[subsection]
\newtheorem{teo}{Teorema}[subsection]

\newtheorem{prop}[teo]{Proposição}
\newenvironment{dem}[1][Demonstra\c{c}{\~a}o]{\textbf{#1:}\ }
{\hfill\rule{1ex}{1ex}}

\begin{document}
\title{Semântica relacional para a lógica proposicional do plausível}
\author{Tiago Augusto dos Santos Boza  \thanks{ Email: boza.tiago@gmail.com. Pós-Graduação em Filosofia, UNESP, FFC - Marília} \\
Hércules de Araujo Feitosa \thanks{Email: haf@fc.unesp.br. Departamento de Matemática, UNESP, FC - Bauru}}

\maketitle \abstract{Esse artigo tem como objetivo apresentar uma semântica de
vizinhança para uma lógica modal, a saber, a lógica proposicional
do plausível. Por ser uma lógica de caráter modal subnormal, não
admite uma semântica de Kripke. Contudo, os autores apresentam neste que as
semânticas de vizinhança, alcançam uma semântica relacional para
a lógica proposicional do plausível.

\medskip

\noindent{\bf Palavras Chave:} lógica modal, lógica proposicional do plausível, semânticas de vizinhança}


\section*{Introdução}

${}$ \hspace {0,5 cm} Não se tem muito claro quais são os objetos específicos de estudo da lógica. Mas, alguns autores ainda tentam explicitar,pelo menos, quais foram e estão sendo estudados.

Desse modo, segundo $\cite{sca76}$, os temas lógicos podem agrupar-se da seguinte maneira:

(i) problemas relacionados com o estudo das inferências válidas e a análise dos conceitos de demonstração e definição;

(ii) problemas que hoje chamamos de semânticos, ligados à análise dos conceitos de significado e verdade;

(iii) análise dos paradoxos lógicos que se apresentam em diferentes ambientes teóricos;

(iv) estudo de alguns conceitos, que segundo $\cite{sca76}$, podem ser denominados de críticos. Estes, por sua vez, visam a tratar dos conceitos de quantidade, números, infinito e etc. 

Estas ideias foram apresentadas em 1976, por $\cite{sca76}$, e por outros autores. Porém, nos parece uma caracterização ainda atual.

Sabemos que, tradicionalmente os sistemas lógicos seguem em suma as assertivas (i) e (ii). Pois, ao se conceber um sistema lógico é necessário que exista uma determinada sintática e uma determinada interpretação semântica.

Agora, cabe-nos aqui a seguinte pergunta: o tratamento dado pela lógica proposicional clássica, com seus conectivos, não possuem um caráter simplista do raciocínio humano? 

Neste caso, nos referimos a um caráter simplista no sentido em que obviamente existem outras formas de raciocínio além dessas?

Esta pergunta e, claramente algumas outras, que motivaram a ideia das chamadas lógicas modais. Estas, por sua vez, inserem um novo operador não definido à partir dos operadores clássicos.

De modo a entendermos quais foram as motivações para a formalização de tais lógicas, apresentamos neste artigo um particular sistema modal, muito difundido na literatura, a saber, o sistema $S_5$. Este
apresenta dois novos operadores $\square$, $\textit{necessário}$, e dualmente $\lozenge$, $\textit{possível}$, já que mostraremos que $\lozenge A \leftrightarrow \neg \square \neg A$. 

E mais, baseando-nos na tradição dos estudos dos sistemas lógicos, trataremos de uma de suas respectivas semânticas, a saber, semântica de Kripke ou semântica dos mundos possíveis. E, portanto, a adequação e/ou equivalência entre essas duas abordagens.

Apresentaremos também uma outra lógica modal, a lógica proposicional do plausível. E, novamente, apresentamos uma interpretação semântica, a partir do contexto algébrico, e demonstramos sua adequação. 

Agora, o que pode nos passar é o seguinte: não há, para esta lógica particular, uma semântica diferente desta apresentada no contexto algébrico e que seja, por exemplo, similar com as semânticas de Kripke?

É exatamente esta a motivação do presente artigo. Portanto, como elemento original deste, apresentamos uma semântica relacional, para a lógica proposicional do plausível. Para tanto, faremos esta apresentação, usando os operadores modais do sistema $S_5$, mostrando sua equivalência dedutiva com a lógica proposicional do plausível e, ao final demonstramos a equivalência entre os dois sistemas, axiomático
e relacional.


\section{Sistema $S_5$ e os Mundos Possíveis}

${}$ \hspace {0,5 cm} Nessa seção, apresentaremos primeiro algumas formalizações básicas
das ideias dos mundos possíveis e, após isso, mostraremos uma construção
para o conhecido sistema $S_5$.

Sabemos que o sistema $S_5$ é determinado semanticamente pelas ideias
de Leibniz de necessidade e possibilidade, dadas por:

\begin{defn}
Uma proposição é necessária se é válida em todos os mundos
possíveis, e é possível se é válida para algum mundo possível.
\end{defn}

Mas ainda não vimos o que são mundos possíveis. Assim, nos próximos
subitens, desenvolveremos algumas ideias semânticas, por meio de uma
definição de verdade em um mundo possível, num modelo para uma linguagem
de necessidade e possibilidade. Isto é, apresentaremos semanticamente
o sistema $S_5$ em termos da noção de verdade e, após isso, sintaticamente por meio de
um sistema dedutivo.

Mas antes, é necessário que apresentemos a linguagem
da qual trataremos a seguir.

\begin{defn}
O conjunto de fórmulas, $\textbf{For($S_5$)}$, é dado por uma extensão
conservativa da definição de fórmulas para a lógica proposicional clássica com o acréscimo das seguintes cláusulas: \\
(i)Se A é uma fórmula, então $\lozenge A$ e $\square A$ também são fórmulas; \\
(ii) O conjunto das fórmulas para $S_5$ é gerado apenas por essas cláusulas.
\end{defn}


\subsection{Verdade e os Mundos Possíveis}

${}$ \hspace {0,5 cm} De acordo com a ideia de Leibniz, mencionada acima, podemos dizer
linguisticamente o seguinte.

\begin{defn}
Uma sentença da forma $\square A$, necessariamente A, é verdadeira
se, e somente se, A é verdadeira em todo mundo possível. A sentença $\lozenge A$,possivelmente A, é verdadeira se, e somente se, A é verdadeira em algum mundo possível.
\end{defn}

Claramente, esta é uma definição que traduz a mesma ideia de Leibniz, apresentada anteriormente, mas devemos notar que, nesse caso, estamos já considerando as sentenças $\square A$ e $\lozenge A$.

Ambas as definições nos passam a ideia de uma coleção de mundos
possíveis, em que sentenças de uma determinada linguagem podem ser verdadeiras
ou falsas, dependendo do mundo em que estejam. 

Desse modo, nosso objetivo é construir um modelo que capture tais ideias.
Faremos tal construção considerando o seguinte: \\

(i) Seja $(P_0, P_1, P_2, ...)$ uma sequência infinita de conjuntos de mundos
possíveis, tais que a condição abaixo esteja satisfeita:
- Para cada número natural $n$, do conjunto de mundos possíveis $P_n$,
destacamos apenas os mundos possíveis em que a variável atômica $p_n$, é
verdadeira.

Mais claramente, segundo a ideia acima, a sequência infinita de conjuntos
de mundos possíveis interpreta cada variável atômica ao estipular em
quais mundos possíveis aquela variável é verdadeira, ou seja, a variável atômica
$p_n$ é verdadeira num mundo possível $\alpha$ se, e somente se, $\alpha$ está no conjunto $P_n$.
Colocando tais ideias numa linguagem formal, segue a seguinte definição.

\begin{defn} Um modelo é um par $<W, P>$, em que: \\
(i) W é um conjunto de mundos possíveis; \\
(ii) P é uma abreviação da sequência infinita $(P_0, P_1, P_2, ...)$ de subconjuntos
de W.
\end{defn}

\begin{defn}
Seja A uma fórmula de $\textbf{For($S_5$)}$ e $\alpha$ é um mundo possível do
modelo $\mathcal{M} = <W, P>$. A fórmula A é válida em $\mathcal{M}$, quando A é verdadeiro no
mundo $\alpha$ de W. Notação: $\models^{\mathcal{M}}_{\alpha}$ A.
\end{defn}

Desse modo, as condições para que as fórmulas sejam válidas são enunciadas
da seguinte forma: \\

(1) $\models^{\mathcal{M}}_{\alpha} p_n$ see $\alpha \in P_n$, para $n = 0, 1, 2, ...$; \\

(2) $\models^{\mathcal{M}}_{\alpha} \top$; \\

(3) $\nvDash^{\mathcal{M}}_{\alpha} \bot$; \\

(4) $\models^{\mathcal{M}}_{\alpha} \neg A$ see $\nvDash^{\mathcal{M}}_{\alpha} A$; \\

(5) $\models^{\mathcal{M}}_{\alpha} A \land B$ see $\models^{\mathcal{M}}_{\alpha} A$ e $\models^{\mathcal{M}}_{\alpha} B$; \\

(6) $\models^{\mathcal{M}}_{\alpha} A \lor B$ see $\models^{\mathcal{M}}_{\alpha} A$ ou $\models^{\mathcal{M}}_{\alpha} B$ ou ambos; \\

(7) $\models^{\mathcal{M}}_{\alpha} A \rightarrow B$ see se $\models^{\mathcal{M}}_{\alpha} A$, então $\models^{\mathcal{M}}_{\alpha} B$; \\

(8) $\models^{\mathcal{M}}_{\alpha} A \leftrightarrow B$ see $\models^{\mathcal{M}}_{\alpha} A$ see $\models^{\mathcal{M}}_{\alpha} B$; \\

(9) $\models^{\mathcal{M}}_{\alpha} \square A$ see $\forall \beta \in \mathcal{M}, \models^{\mathcal{M}}_{\beta} A$; \\

(10) $\models^{\mathcal{M}}_{\alpha} \lozenge A$ see $\exists \beta \in \mathcal{M}, \models^{\mathcal{M}}_{\beta} A$. \\

A partir de agora, tomaremos a seguinte definição para simplificar as notações.

\begin{defn}
Dizemos que A é válida, o que é denotado por $\models A$, se para
todo modelo $\mathcal{M}$ e qualquer mundo possível $\alpha$ de $\mathcal{M}$, temos que $\models^{\mathcal{M}}_{\alpha} A$.
\end{defn}

Postas essas considerações, abaixo apenas apresentaremos alguns resultados sabidamente válidos. A cada resultado corresponde uma determinada denominação que será apresentada precedendo o esquema de fórmula.

\begin{prop} 
O esquema $\textbf{T}: \square A \rightarrow A$ é válido em $\mathcal{M}$.
\end{prop}

\begin{prop}
O esquema $\textbf{5}: \lozenge A \rightarrow \square \lozenge A$ é válido em $\mathcal{M}$.
\end{prop}

\begin{prop}
O esquema $\textbf{K}: \square(A \rightarrow B) \rightarrow (\square A \rightarrow \square B$ é válido em $\mathcal{M}$.
\end{prop}

O próximo resultado se refere à $\textit{regra de necessitação}$.[

\begin{prop}
Se $\models A$, então $\models \square A$.
\end{prop}

Já, o próximo item, refere-se à uma definição do operador $\lozenge$ em termos dos operadores $\square$ e $\neg$, denotada por $\textbf{Df} \lozenge: \lozenge A \leftrightarrow \neg \square \neg A$. Mais especificamente, esta corresponde à definição dual do operador $\lozenge$.

\begin{prop}
O esquema $\textbf{Df} \lozenge: \lozenge A \leftrightarrow \neg \square \neg A$ é válido em $\mathcal{M}$.
\end{prop}

Abaixo, alguns resultados muito importantes. O primeiro relaciona as tautologias com os modelos $\mathcal{M}$. E, o segundo garante a validade da regra de $Modus Ponens$.

\begin{prop}
Se A é uma tautologia, então $\models A$. 
\end{prop}

\begin{prop}
Se $\models A \rightarrow B$ e $\models A$, então $\models B$.
\end{prop}


\subsection{Sistema Axiomático para $S_5$}

${}$ \hspace {0,5 cm} Nessa seção apresentaremos uma versão axiomática do Sistema $S_5$.

\begin{defn}
o Sistema $S_5$ é uma extensão conservativa da lógica proposicional
clássica pelo acréscimo dos seguintes itens:\\
- O alfabeto é obtido pelo acréscimo dos operadores $\square$ e $\lozenge$, não definidos à partir dos outros operadores; \\
- Os axiomas de $S_5$ é determinado pelos seguintes esquemas de axiomas: \\

(i) Axiomas da lógica proposicional clássica; \\

(ii) $\textbf{T}: \square A \rightarrow A$; \\

(iii) $\textbf{5}: \lozenge A \rightarrow \square \lozenge A$; \\

(iv) $\textbf{K}: \square (A \rightarrow B) \rightarrow (\square A \rightarrow \square B)$; \\

(v) $\textbf{Df}\lozenge: \lozenge A \leftrightarrow \neg \square \neg A$. \\

E pelas seguintes regras de dedução: \\

(vi) $\textbf{RN}: \vdash A / \vdash \square A$; \\

(vii) $\textbf{MP}: A \rightarrow B, A / B$.
\end{defn}

Além disso, definimos o seguinte neste sistema.

\begin{defn}
Um teorema é uma sentença que pode ser obtida a partir
dos axiomas e das regras de inferência.
\end{defn}

Abaixo, apenas apresentaremos alguns resultados que são sabidamente válidos e que também podem ser encarados como regras de dedução no sistema $S_5$.

\begin{prop}
Para $n \geq 0$, vale $\textbf{RPL}: A_1 \rightarrow A_2 \rightarrow A_3 \rightarrow ... \rightarrow A_n / A$.
\end{prop}

\begin{prop}
O esquema $\textbf{T}\lozenge: A \rightarrow \lozenge A$ é um teorema de $S_5$.
\end{prop}

\begin{prop}
O esquema $\textbf{D}: \square A \rightarrow \lozenge A$ é um teorema de $S_5$.
\end{prop}

\begin{prop}
O esquema $\textbf{B}: A \rightarrow \square \lozenge A$ é um teorema de $S_5$.
\end{prop}

\begin{prop}
A regra $\textbf{RM}: \vdash A \rightarrow B / \vdash \square A \rightarrow \square B$ pode ser obtida em $S_5$.
\end{prop}

\begin{prop}
A regra $\textbf{RE}: \vdash A \leftrightarrow B / \vdash \square A \leftrightarrow \square B$ pode ser obtida em $S_5$.
\end{prop}

E por fim, uma definição dual de $\square$, à partir de $lozenge$ e $neg$.

\begin{prop}
O esquema $\textbf{Df}\square: \square A \leftrightarrow \neg \lozenge \neg A$ é um teorema de $S_5$.
\end{prop}

E mais alguns resultados.

\begin{prop}
O esquema $\textbf{B}\lozenge: \lozenge \square A \rightarrow A$ é um teorema de $S_5$.
\end{prop}

\begin{prop}
O esquema $\textbf{5}\lozenge: \lozenge \square A \rightarrow \square A$ é um teorema de $S_5$.
\end{prop}

\begin{prop}
O esquema $\textbf{4}: \square A \rightarrow \square \square A$ é um teorema de $S_5$.
\end{prop}

\begin{prop}
O esquema $\textbf{4}\lozenge: \lozenge \lozenge A \rightarrow \lozenge A$ é um teorema de $S_5$.
\end{prop}

\begin{prop}
Para todo $n \geq 0$ a regra $\textbf{RK}: (A_1 \land A_2 \land A_3 \land ... \land A_n) \rightarrow A / (\square A_1 \land ... \land \square A_n) \rightarrow \square A$ pode ser obtida em $S_5$.
\end{prop}

Notemos que um caso particular de $\textbf{RK}$ recebe uma denominação diferenciada,
o caso em que $n = 0$, claramente é a regra $\textbf{RN}$. Um outro caso
particular de $\textbf{RK}$ que recebe uma denominação diferenciada é o caso em que
$n = 2$ e a denominação é $\textbf{RR}$. E também alguns resultados que advém imediatamente
destes são: \\

(i) $\textbf{N}: \square \top$, resultado imediato de $\textbf{RN}$;\\

(ii) $\textbf{M}. \square (A \land B) \rightarrow (\square  A \land \square B)$, imediatamente de $\textbf{RM}$;\\

(iii) $\textbf{C}. (\square A \land \square B \rightarrow \square (A \land B)$, provém imediatamente de $\textbf{RR}$. \\


\subsection{Adequação para $S_5$}

${}$ \hspace {0,5 cm} Nesta seção mostraremos uma adequação entre o sistema apresentado
acima, $S_5$, e as semânticas de Kripke.

\begin{defn}
Dizemos que um conjunto de fórmulas $\Gamma$ da linguagem $\mathbb{L}(\neg, \land)$ da lógica proposicional clássica é dito completo ou maximalmente consistente quando para toda fórmula
A de $\mathbb{L}$, ou ocorre A ou ocorre $\neg A$ em $\Gamma$.
\end{defn}

\begin{defn}
Dada uma lógica modal $\mathfrak{L}$, dizemos que um conjunto de fórmulas $\Gamma$ é $\mathfrak{L}_{\land}$ inconsistente se para $B_1, ..., B_n \in \Gamma$ e $\vdash_{\mathfrak{L}} \neg (B_1 \land ... \land B_n)$. Caso contrário, $\Gamma$ é dito $\mathfrak{L}_{\land}$ consistente.
\end{defn}

\begin{defn}
Um modelo de Kripke para uma lógica modal $\mathfrak{L}$ é uma terna da seguinte forma $\mathcal{M} = <W, R, e>$, em que: \\

- W é o conjunto de todos os mundos possíveis;\\

- R é uma relação binária, denominada de relação de acessibilidade, entre elementos de W; \\

- e é uma valoração tal que e: $W \rightarrow \mathrm{P}(Var(\mathfrak{L}))$. \\

Quando uma fórmula A de $\mathfrak{L}$ é válida num mundo $w \in W$, denotaremos por $\vDash^{\mathcal{M}}_{w} A$.
\end{defn}

\begin{defn}
Uma fórmula de $\mathfrak{L}$ é válida se a mesma é válida em todo mundo w de $\mathcal{M}$.
\end{defn}

Dessa maneira, apresentamos uma interpretação semântica para o sistema S5.
Esta semântica é constituída por elementos que chamamos de mundos possíveis,
por uma relação binária definida entre esses elementos e uma valoração.
Esta semântica é conhecida como Semântica de Kripke.

Para mostrarmos o que queremos adotaremos, neste trabalho, o que chamamos de modelo
canônico para a referida semântica, um caso particular. Mais especificamente,
consideraremos uma estrutura um pouco diferente daquela apresentada
acima.

\begin{defn}
Um modelo canônico para uma lógica modal $\mathfrak{L}$ é uma terna da seguinte forma $\mathcal{M} = <W_L, R_L, e_L>$, em que: \\

- $W_L$ é o conjunto de todos os conjuntos $\Gamma$ de fórmulas, que são L completos;\\

- $R_L$ é uma relação tal que: dados $\Gamma$, $\Delta$ um conjunto de fórmulas de $W_L$:
$\Gamma R_L \Delta$ see para toda fórmula A, se $\square A \in \Gamma$, então $A \in \Delta$; \\

- $e_L$ é uma valoração tal que: $e_L(\Gamma) = \lbrace p : p \in \Gamma \rbrace$.
\end{defn}

Agora, as condições para que as fórmulas sejam válidas nessa estrutura
são enunciadas da seguinte forma: \\

(i) As condições de (1) a (8) enunciadas nas seções anteriores do presente artigo,
acrescidas das seguintes condições: \\

$(9)^{*}$ $\models_{\mathcal{M}} \square A$ see $\forall z \in \mathcal{M}$ tal que $w R_L z$, $\models^{\mathcal{M}}_{z} A$; \\

$(10)^{*}$ $\models_\mathcal{M} \lozenge A$ see $\exists z \in \mathcal{M}$ tal que $w R_L z$, $\models^{\mathcal{M}}_{z} A$. \\

\begin{defn}
Dizemos que uma classe de frames $\mathcal{C}$ é correta para uma lógica
modal $\mathfrak{L}$ quando: $\vdash_{\mathfrak{L}} A$ $\Rightarrow$ $\forall \mathcal{M} \in \mathcal{C}$, $\mathcal{M} \vDash A$.
\end{defn}

\begin{defn}
Uma classe de frames $\mathcal{C}$ é completa para uma lógica modal $\mathfrak{L}$ quando: $\forall \mathcal{M} \in \mathcal{C}$, $\mathcal{M} \vDash A$ $\Rightarrow$ $\vdash_{\mathfrak{L}} A$.
\end{defn}

\begin{defn}
Dizemos que uma relação R é uma relação de equivalência
quando as duas condições abaixo são satisfeitas:\\

(i) R é reflexiva, i.e., para todo $\alpha$, $\alpha R \alpha$; \\

(ii) R é euclidiana, i.e., para todo $alpha$, $\beta$ e $\gamma$, se $\alpha R \beta$ e $\alpha R \gamma$, então $\beta R \gamma$.
\end{defn}

Segundo $\cite{eps95}$ a relação dada num modelo para $S_5$ é reflexiva e euclidiana, respectivamente, quando: \\

(i) R é reflexiva see $\models^{\mathcal{M}}_{\alpha} A \square A \rightarrow A$; \\

(ii) R é euclidiana see $\models^{\mathcal{M}}_{\alpha} \lozenge A \rightarrow \square \lozenge A$. \\

Desse modo, basta mostrarmos estes itens.

\begin{prop}
São válidas as seguintes sentenças: \\

(i) $\models^{\mathcal{M}}_{\alpha} A \square A \rightarrow A$; \\

(ii) $\models^{\mathcal{M}}_{\alpha} A \square A \rightarrow A$.
\end{prop}

\begin{prop}
A classe de frames que tem uma relação de equivalência são corretas para $S_5$.
\end{prop}

\begin{prop}
Se $\Gamma$ é $\mathfrak{L}_{\land}$ consistente, então $\Gamma \cup \lbrace \neg A \rbrace$ ou $\Gamma \cup \lbrace A \rbrace$ é $\mathfrak{L}_{\land}$ consistente. Assim, se $\Gamma$ é $\mathfrak{L}_{\land}$ consistente e completo, então para cada fórmula A, ou A ou $\neg A$ está em $\Gamma$.
\end{prop}

\begin{prop}
Se $\Omega$ é $\mathfrak{L}_{\land}$ consistente, então existe $\Gamma$ tal que, $\Omega \subseteq \Gamma$ e $\Gamma$ é $\mathfrak{L}_{\land}$ consistente e completo.
\end{prop}

\begin{prop}
Se $\Gamma$ é $\mathfrak{L}_{\land}$ consistente e completo, então: \\

(i) Teo($\mathcal{L}) \subseteq \Gamma$ e $\Gamma$ é fechado para a regra modus ponens; \\

(ii) Se $A \vdash_{LPC} B$, então: se $\square A \in \Gamma$, então $\square B \in \Gamma$ e se $\neg \square B \in \Gamma$, então $\neg \square A \in \Gamma$. \\

Notar que Teo$(\mathcal{L})$ são os teoremas da lógica proposicional clássica.
\end{prop}

Voltemos aos modelos canônicos.

\begin{prop}
Seja $\Gamma \in W_L$ tal que para todo $\Omega \in W_L$, com $\Gamma R_L \Omega$, implica $B \in \Omega$. Então, $\square B \in \Gamma$.
\end{prop}

\begin{prop}
As duas sentenças abaixo são verdadeiras: \\

(i) $<W_L, R_L, e_L, \Gamma>$ $\models A$ $\Leftrightarrow$ $A \in \Gamma$; \\

(ii) $\vdash_{\mathcal{L}} A$ $\Leftrightarrow$ $<W_L, R_L, e_L>$ $\models$ $A$. \\
\end{prop}

Desses resultados segue o queremos.

\begin{prop} $\textit{Adequação das Semânticas de Kripke}$ As classes de frames que possuem uma relação de equivalência são modelos adequados para o sistema $S_5$.
\end{prop}


\section{A motivação quantificacional}

${}$\hspace{0,5cm} A formalização lógica quantificacional dos espaços pseudo-topológicos surgiu em $\cite{gra99}$ como uma extensão da lógica clássica de primeira ordem com igualdade $\mathcal{L}$. Para detalhes sobre $\mathcal{L}$ ver $\cite{eft84}$ ou $\cite{fep05}$. Foi denominada por Grácio de lógica do plausível e denotada por $\mathcal{L}(P)$. \\

Dada $\mathcal{L}$, a lógica estendida $\mathcal{L}(P)$ é determinada pelos seguintes acréscimos: \\

- a linguagem $L$ de $\mathcal{L}(P)$ conta com um novo símbolo de quantificador $P$ e sentenças do tipo $Px \varphi(x)$ são bem formadas em  $\mathcal{L}(P)$. \\

- axiomas específicos do quantificador $P$:

$(A_1)$ $Px \varphi(x) \land Px \psi(x) \to Px(\varphi(x) \land Px \psi(x))$

$(A_2)$ $Px \varphi(x) \land Px \psi(x) \to Px(\varphi(x) \lor Px \psi(x))$

$(A_3)$  $\forall x \varphi(x) \rightarrow Px \varphi(x)$

$(A_4)$ $Px \varphi(x) \rightarrow \exists x \varphi(x)$

$(A_5)$ $\forall x (\varphi(x) \to \psi(x)) \to (Px \varphi(x) \to Px \psi(x))$

$(A_6)$ $Px \varphi(x) \to Py \varphi(y)$, se $y$ é livre para $x$ em $\varphi(x)$. \\

- regradas de dedução:

(MP) \emph{Modus Ponens}: $\varphi, \varphi \to \psi \vdash \psi$

(Gen) Generalização: $\varphi \vdash \forall x \varphi(x)$. \\

Os demais conceitos sintáticos usuais como sentença, demonstração, teorema, dedução, consistência e outros são definidos do modo padrão. \\

As estruturas adequadas para $\mathcal{L}(P)$, denominadas por Grácio de estrutura do plausível, também são extensões das estruturas de primeira ordem $\mathcal{A}$. \\

Assim dada uma estrutura $\mathcal{A}$, consideremos que o seu domínio é denotado por $A$. Uma estrutura do plausível, denotada por $\mathcal{A}^{\Omega}$, é determinada a partir $\mathcal{A}$ pelo acréscimo de um espaço pseudo-topológico $\Omega$ sobre o universo $A$.

A interpretação dos símbolos de relação, função e constante é a mesma de $\mathcal{L}$ com relação à $\mathcal{A}$.

\begin{defn} A satisfação de uma sentença do tipo $Px \varphi(x)$ em $\mathcal{A}^{\Omega}$ e definida indutivamente por:

- se $\varphi$ é uma fórmula cujas variáveis livres estão em $\{x\} \cup \{y_1, ..., y_n\}$ e $\overline{a} =
(a_1, ..., a_n)$ é uma sequência de elementos de $A$, então: 
$$\mathcal{A}^{\Omega} \vDash Px \varphi[x, \overline{a}] \Leftrightarrow \{b \in A : \mathcal{A}^{\Omega} \vDash
[b, \overline{a}]\} \in \Omega.$$
\end{defn} 

Da maneira usual, para a sentença $Px \varphi(x)$:
$$\mathcal{A}^{\Omega} \vDash Px \varphi(x) \Leftrightarrow \{a \in A : \mathcal{A}^{\Omega} \vDash \varphi(a)\} \in \Omega.$$

As outras noções semânticas como modelo, validade, implicação lógica, entre outras, são apropriadamente adaptadas a partir da interpretação de $\mathcal{L}$ em $\mathcal{A}$. \\

Grácio $\cite{gra99}$ provou que as estruturas do plausível são modelos corretos e completos para $\mathcal{L}(P)$.


\subsection{Uma formalização axiomática e proposicional}

${}$\hspace{0,5cm} Apresentamos a versão proposicional da lógica do plausível, conforme $\cite{fng09}$, que procura formalizar os aspectos de um espaço pseudo-topológico no contexto lógico proposicional, sem a presença dos quantificadores, inclusive o do plausível, que é substituído por um operador unário. \\

Esta lógica será denotada por $\mathbb{L}(\nabla)$. Como mencionado, ela estende a lógica proposicional clássica (LPC) na linguagem $L(\neg, \land, \lor, \to)$ com o acréscimo do operador unário $\nabla$, donde obtermos a linguagem proposicional $L(\neg, \land, \lor, \to, \nabla)$.

A lógica fica determinada pelo seguinte: \\

- Axiomas:

$(Ax_0)$ LPC

$(Ax_1)$ $(\nabla \varphi \land \nabla \psi) \to \nabla(\varphi \land \psi)$

$(Ax_2)$ $\nabla(\varphi \lor \neg \varphi)$ 

$(Ax_3)$ $\nabla \varphi \to \varphi$. \\

- Regras de dedução:

$(MP)$ \emph{Modus Ponens}

$(R\nabla)$ $\vdash \varphi \to \psi\ / \vdash \nabla \varphi \to \nabla \psi$. \\

A intuição para este operador de plausível é de algo que pode ser explicado numa teoria sem provocar inconsistência. Assim, o plausível não se assemelha ao possível, pois podemos ter: ``é possível que chova amanhã'' e ``é possível que não chova amanhã'', porém ``não é possível que chova e não chova amanhã''.

O conceito está vinculado, por exemplo, com a existência de uma demonstração de um fato numa teoria consistente. Assim, não pode haver uma prova de $\varphi$ e uma outra de $\neg \varphi$.

Diante disso, $\varphi$ e $\psi$ são plausíveis se, e somente se, $\varphi \land \psi$ é plausível. Toda tautologia é plausível e se $\varphi$ é plausível, então vale a proposição $\varphi$. A regra $(R\nabla)$ diz que se há uma prova de $\varphi \to \psi$ e $\varphi$ é plausível, então também $\psi$ é plausível. \\

	Pode ser demonstrado o seguinte.
	
\begin{prop} (i) $\vdash \neg \nabla \bot$

(ii) $\vdash \nabla \varphi \to \nabla (\varphi \lor \psi)$

(iii) $\vdash \varphi \Rightarrow\ \vdash \nabla \varphi$

(iv) $\vdash \varphi \to \neg \nabla \neg \varphi$

(v) $\vdash \nabla \varphi \to \neg \nabla \neg \varphi$

(vi) $\vdash \nabla \neg \varphi \to \neg \nabla \varphi$.
\end{prop}

\begin{prop} 
$\nabla \varphi \to \nabla(\varphi \lor \psi) \Leftrightarrow (\nabla \varphi \lor \nabla \psi) \to \nabla(\varphi \lor \psi)$. \\
\begin{dem} $(\Rightarrow)$ Da hipótese, $\nabla \varphi \to \nabla(\varphi \lor \psi)$ e $\nabla \psi \to \nabla(\varphi \lor \psi)$. Daí, $(\nabla \varphi \lor \nabla \psi) \to \nabla(\varphi \lor \psi)$.
$(\Leftarrow)$ Como $\nabla \varphi \to (\nabla \varphi \lor \nabla \psi)$, então segue da hipótese que $\nabla \varphi \to \nabla(\varphi \lor \psi)$. 
\end{dem}
\end{prop}


\subsection{Álgebras dos espaços pseudo-topológicos}

${}$\hspace{0,5cm} Agora, apresentamos a álgebra do plausível, que correspondem à versão algébrica da lógica da seção anterior, conforme \cite{fng09}.

\begin{defn} Álgebra do Plausível é uma estrutura $\mathbb{P} = (P, 0, 1, \land, \lor, \sim, \sharp)$, em que $(P, 
0, 1, \land, \lor, \sim)$ é uma álgebra de Boole e $\sharp$ é o operador do plausível, sujeito a: 

$(a_1)$  $\sharp a \land \sharp b \leq \sharp (a \land b)$

$(a_2)$  $\sharp a \leq \sharp (a \lor b)$

$(a_3)$  $\sharp a \leq a$

$(a_4)$  $\sharp 1 = 1$. 	\end{defn}

\begin{defn} Um elemento $a \in P$ é plausível $a \neq 0$ e $\sharp a = a$.	\end{defn}  

Embora $\sharp 0 = 0$, por definição, $0$ não é plausível.

\begin{prop} Se $\mathbb{P} = (P, 0, 1, \land, \lor, \sim, \sharp)$ é uma álgebra do plausível e $a, b \in P$, então:

(i) $\sharp a \leq \sharp(a \lor b)$

(ii) $a \leq b \Rightarrow \sharp a \leq \sharp b$

(iii) $\sharp a \lor \sharp b \leq \sharp (a \lor b)$. 
\end{prop}

\begin{prop} Para cada álgebra do plausível $\mathbb{P} = (P, 0, 1, \land, \lor, \sim, \sharp)$ existe um monomorfismo $h$  de $P$ num espaço pseudo-topológico de conjuntos definidos em $\mathcal{P}(\mathcal{P}(P))$.
\end{prop}

Em $\cite{fng09}$ há uma demonstração da adequação de $\mathbb{L}(\nabla)$ com relação às álgebras do plausível $\mathbb{P}$. \\

A interpretação de $\mathbb{L}(\nabla)$ numa pseudo-topologia $(E, \Omega)$ é uma função $v$ com as seguintes características:
	
As sentenças universais que valem para todos os indivíduos, são agora tomadas pelas tautologias, que são interpretados no universo $E$ da pseudo-topologia. 

Se $\varphi$ e $\psi$ estão $\Omega$ e, portanto, são ubíquos, então o mesmo vale  $\varphi \land \psi$.

Se $\varphi$ ou $\psi$ está $\Omega$, então o mesmo vale  $\varphi \lor \psi$. Esta condição implica que se $\varphi, \psi \in \Omega$, então $\varphi \lor \psi \in \Omega$, mas é um pouco mais básica. Aqui está uma variação entre os dois sistemas.

Se $\varphi$ e $\psi$ são equivalentes, então uma delas está numa pseudo-topologia se, e somente se, também a outra está.

A interpretação de $\bot$ não pode estar $\Omega$. Para tanto, usamos o axioma $(Ax_3)$ $\nabla \bot \to \bot$ que implica esta condição no contexto proposicional. Mas também tem mais exigência que $(E_4)$. \\

Talvez possamos, em algum momento, refinar estas variações, mas por ora assumamos a lógica $\mathbb{L}(\nabla)$.


\section{Semântica Relacional para a Lógica Proposicional do Plausível}

${}$ \hspace {0,5 cm} Nessa parte do presente artigo apresentaremos uma semântica relacional para
a lógica proposicional do plausível, e também apresentaremos a adequação entre esta estrutura semântica
relacional e a axiomática, introduzida anteriormente.

Para isto, consideremos a seguinte lógica proposicional modal.

\begin{defn}
A lógica proposicional modal do plausível, $\mathbb{LP}(\square)$, é uma lógica
modal que estende conservativamente a lógica proposicional clássica pelo
acréscimo do operador de necessidade, $\square$, com os seguintes axiomas e regras
específicos para este operador: \\

$\textbf{C}: \square A \land \square B \rightarrow \square (A \land B)$ \\

$\textbf{H}: \square A \lor \square B \rightarrow \square (A \lor B)$ \\

$\textbf{T}: \square A \rightarrow A$ \\

$\textbf{N}: \square \top$ \\

E a regra de inferência: \\

$\textbf{RE}: \vdash A \leftrightarrow B / \vdash \square A \leftrightarrow \square B$
\end{defn}

Com esta caracterização mostraremos, agora, que este sistema e a lógica
proposicional do plausível são dedutivamente equivalentes.

\begin{prop}
Os sistemas $\mathbb{LP}(\square)$ e $\mathbb{L}(\nabla)$ são dedutivamente equivalentes.
\begin{dem}
Por abstração é imediato que os esquemas de axiomas $(Ax_1)$ e $\textbf{C}$, $(Ax_3)$ e $\textbf{T}$, $(Ax_4)$ e $\textbf{N}$, são equivalentes, respectivamente. Já para $\textbf{H}$, sabemos
que $\nabla A \rightarrow \nabla (A \lor B)$ é equivalente à $\nabla A \lor \nabla B \rightarrow \nabla (A \lor B)$ e, assim, obtemos que os esquemas de axiomas são equivalentes. Desse modo, falta-nos mostrar que as regras de inferência $(R \nabla)$ e $\textbf{RE}$ são equivalentes. Este resultado segue, também imediatamente, por abstração. Portanto, os sistemas $\mathbb{LP}(\square)$ e $\mathbb{L}(\nabla)$ são dedutivamente equivalentes.
\end{dem}
\end{prop}

A seguir apresentaremos um modelo que será uma frame para $\mathbb{L}(\nabla)$, via $\mathbb{LP}(\square)$.

\begin{defn} Um modelo relacional $\mathcal{M}$ para a lógica proposicional do plausível é uma terna $\mathcal{M} = <W, S, V>$, em que:\\

- W é um conjunto não-vazio de mundos possíveis;\\

- S uma função que associa cada $\alpha \in W$ um conjunto de subconjuntos de W, i.e., $S(\alpha) \subseteq \mathbb{P}(U)$, de tal modo que satisfaz as seguintes condições:\\

(c) $X \in S(\alpha)$ e $Y \in S(\alpha)$ $\Rightarrow$ $X \cap Y \in S(\alpha)$; \\

(h) $X \in S(\alpha)$ ou $Y \in S(\alpha)$ $\Rightarrow$ $X \cup Y \in S(\alpha)$; \\

(t) $X \in S(\alpha)$ $\Rightarrow$ $\alpha \in X$; \\

(n) $W \in S(\alpha)$; \\

- V é uma valoração em U, i.e., uma função do conjunto das fórmulas atômicas em $\mathbf{P}(U)$.
\end{defn}

Agora, temos de definir condições para que as fórmulas sejam válidas nesta nova estrutura.

\begin{defn} Seja $\mathcal{M}$ um modelo e $\alpha \in W$. Uma fórmula A é verdadeira
no mundo $\alpha$, o que será denotado por $\models^{\mathcal{M}}_{\alpha} A$, quando: \\

- as condições para os operadores clássicos seguem como definido anteriormente para o sistema $S_5$;\\ 

E acrescentamos o seguinte:\\

- $\models^{\mathcal{M}}_{\alpha} p_i$ $\Leftrightarrow$ $\alpha \in V(p_i)$,  se $p_i$ é uma variável proposicional; \\

- $\models^{\mathcal{M}}_{\alpha} \square A$ $\Leftrightarrow$ $\parallel A \parallel ^{\mathcal{M}} \in S(\alpha)$, em que $\parallel A \parallel ^{\mathcal{M}} = \lbrace \alpha \in W:$ $\models^{\mathcal{M}}_{\alpha} A \rbrace$; \\

O conjunto $\parallel A \parallel ^{\mathcal{M}}$ será denominado por conjunto verdade de A em $\mathcal{M}$.
\end{defn}

Em alguns momentos, quando não causar estranheza, denotaremos o conjunto verdade apenas por $\parallel A \parallel$.

\begin{defn}
Uma fórmula A é válida num modelo $\mathcal{M}$, quando é verdadeira para todo $\alpha \in W$. A fórmula A é válida se é verdadeira em todo modelo $\mathcal{M}$.
\end{defn}

Denotaremos que uma fórmula A é válida num modelo $\mathcal{M}$ por $\models^{\mathcal{M}} A$, e que A é válida apenas por $\models A$.

\begin{defn}
Se $\Gamma$ é um conjunto de fórmulas e $\mathcal{M}$ é um modelo, então dizemos que $\models^{\mathcal{M}} \Gamma$ see $\models^{\mathcal{M}} \Gamma$, para cada $A \in \Gamma$.
\end{defn}

\begin{defn}
Seja $\Gamma \cup \lbrace A \rbrace$ um conjunto de fórmulas. Dizemos que $\Gamma$ implica A, ou que A é uma conseqüência global semântica de $\Gamma$, quando para todo modelo $\mathcal{M}$, temos: $\models^{\mathcal{M}} \Gamma$ $\Rightarrow$ $\models^{\mathcal{M}} A$.
\end{defn}

Denotaremos uma conseqüência global semântica por: $\Gamma \models A$.

Agora, precisamos de um resultado a ser demonstrado.

\begin{prop}
Seja $\mathcal{M}$ um $\mathbb{L}(\nabla)$-modelo e A e B fórmulas quaisquer. Então: \\

(i) $\parallel \neg A \parallel$ = $- \parallel A \parallel$; \\

(ii) $\parallel A \land B \parallel$ = $\parallel A \parallel \cap \parallel B \parallel$; \\

(iii) $\parallel A \lor B \parallel$ = $\parallel A \parallel \cup \parallel B \parallel$; \\

(iv) $\parallel A \rightarrow B \parallel$ = $- \parallel A \parallel \cup \parallel B \parallel$; \\

(v) $\parallel A \leftrightarrow B \parallel$ = $(- \parallel A \parallel \cup \parallel B \parallel) \cap (- \parallel B \parallel \cup \parallel A \parallel)$; \\

(vi) $\parallel \square A \parallel$ = $\lbrace \alpha \in W$ : $\parallel A \parallel \in S(\alpha) \rbrace$. \\
Observação: o conjunto, $- \parallel A \parallel$, equivale ao conjunto $W - \parallel A \parallel$, isto é, é o complemento de $\parallel A \parallel$ relativo à W. \\
\begin{dem}
As condições de (i) a (v) seguem imediatamente do caráter Booleano da lógica clássica e da teoria dos conjuntos. Dessa maneira, basta mostrarmos a condição (vi). Por definição, $\alpha \models \square A$ see $\alpha \in$ $\parallel \square A \parallel$. E isto significa, portanto, que $\lbrace \alpha \in W$ : $\parallel A \parallel$ $\in S(\alpha) \rbrace$ = $\parallel \square A \parallel$.
\end{dem}
\end{prop}

E assim, mostraremos abaixo parte do que queremos. A saber, a correção da lógica $\mathbb{LP}(\square)$ relativa à $\mathcal{M}$.

\begin{prop} $\textit{(Correção)}$ $\Gamma \vdash A$ $\Rightarrow$ $\Gamma \models A$. \\
\begin{dem}
Consideremos $\mathcal{M}$ um $\mathbb{L}(\nabla)$-modelo e façamos a demonstração por indução sobre o comprimento da dedução. \\
- Para $n = 1$, temos que A pertence à $\Gamma$ ou A é um axioma. Desse modo, façamos em duas partes:\\

(i) Se $A \in \Gamma$ ou A é um axioma do cálculo proposicional clássico, nada temos a demonstrar; \\

(ii) Agora, se A é um dos axiomas modais, então: \\

(C) Seja $\alpha \in W$, tal que $\models^{\mathcal{M}}_{\alpha} \square A \land \square B$. Daí, por definição, temos que $\models^{\mathcal{M}}_{\alpha} \square A$ e $\models^{\mathcal{M}}_{\alpha} \square B$ e, assim, segue que, $\parallel A \parallel$ $\in S(\alpha)$ e $\parallel B \parallel$ $\in S(\alpha)$. Pela condição (c) obtemos $\parallel A \parallel \cap \parallel B \parallel$ $\in S(\alpha)$ e, por definição, $\parallel A \land B \parallel$ $\in S(\alpha)$. Desse modo, $\models^{\mathcal{M}}_{\alpha} \square (A \land B)$. Portanto, C vale. \\

(H) Seja $\alpha \in W$, tal que $\models^{\mathcal{M}}_{\alpha} \square A \lor \square B$. Daí, por definição, temos que $\models^{\mathcal{M}}_{\alpha} \square A$ ou $\models^{\mathcal{M}}_{\alpha} \square B$ e, assim, segue que, $\parallel A \parallel$ $\in S(\alpha)$ ou $\parallel B \parallel$ $\in S(\alpha)$. Pela condição (h) obtemos $\parallel A \parallel \cup \parallel B \parallel$ $\in S(\alpha)$ e, por definição, $\parallel A \lor B \parallel$ $\in S(\alpha)$. Desse modo, $\models^{\mathcal{M}}_{\alpha} \square (A \lor B)$. Portanto, H vale. \\

(T) Seja $\alpha \in W$, tal que $\models^{\mathcal{M}}_{\alpha} \square A$. Dessa maneira, por definição, $\parallel A \parallel$ $\in S(\alpha)$ e, pela condição (t), obtemos que $\alpha \in$ $\parallel A \parallel$. Como isto equivale à $\models^{\mathcal{M}}_{\alpha} A$ segue, portanto, que T é válido. \\

- Agora, tomemos como hipótese de indução que: $\Gamma \vdash A_n$ $\Rightarrow$ $\Gamma \models A_n$, para todo $n \leq k$. \\

Assim, há três possibilidades para o passo seguinte, $k+1$, da indução:\\
(i) $A_{k+1}$ é uma premissa;\\
(ii) $A{k+1}$ é um esquema de axiomas;\\
(iii) $A{k+1}$ é deduzida a partir das regras MP ou RE; \\

Para os itens (i) e (ii), nada temos a demonstrar, pois estes itens ficam como na base da indução. E mais, para a regra MP, sabemos que a mesma preserva a validade. Desse modo, falta-nos analisar a regra RE.\\
(RE) Consideremos a regra RE e suponhamos que $\vdash A \leftrightarrow B$. Desse modo, pela Hipótese de Indução, segue que $\models A \leftrightarrow B$. Logo, A e B são equivalentes e, obtemos que, $\parallel A \parallel$ = $\parallel B \parallel$. Daí, para todo $\alpha \in W$, $\parallel A \parallel$ $\in S(\alpha)$ see $\parallel B \parallel$ $\in S(\alpha)$. Assim, $\models^{\mathcal{M}}_{\alpha} \square A$  see $\models^{\mathcal{M}}_{\alpha} \square B$. Por definição, temos $\models^{\mathcal{M}}_{\alpha} \square A \leftrightarrow \square B$. Portanto, RE preserva a validade.
\end{dem}
\end{prop}

Agora, daremos algumas definições e mostraremos alguns resultados para, por fim, demonstrarmos a completude de nosso sistema, relativo à semântica $\mathcal{M}$.

\begin{defn}
Um conjunto de fórmulas $\Delta$ é maximalmente consistente quando $\Delta$ é consistente e nenhuma extensão conservativa de $\Delta$ é consistente.
\end{defn}

A demonstração do resultado abaixo pode ser encontrado em $\cite{che99}$.

\begin{prop} $\textit{(Lindenbaum)}$
Todo conjunto maximalmente consistente $\Gamma$ pode ser estendido à um conjunto maximalmente consistente $\Delta$.
\end{prop}

Assim como a demonstração da completude se deu para o $S_5$, aqui passaremos a considerar os modelos canônicos. Consideremos, então, o conjunto de todos os conjuntos $\mathbb{L}(\nabla)$-Maximalmente
Consistente, que será denotado por $\chi$.

\begin{defn}
O conjunto de demonstrações de A é o conjunto $\vert A \vert$ = $\lbrace \Gamma \in \chi$ : $A \in \Gamma \rbrace$.
\end{defn}

\begin{prop}
Sejam A e B fórmulas quaisquer. Então:\\

(i) $\vert \neg A \vert$ = $- \vert A \vert$; \\

(ii) $\vert A \land B \vert$ = $\vert A \vert \cap \vert B \vert$; \\

(iii) $\vert A \lor B \vert$ = $\vert A \vert \cup \vert B \vert$; \\

(iv) $\vert A \rightarrow B \vert$ = $- \vert A \vert \cup \vert B \vert$; \\

(v) $\vert A \leftrightarrow B \vert$ = $(- \vert A \vert \cup \vert B \vert) \cap (- \vert B \vert \cup \vert A \vert)$; \\

(vi) $\vert A \vert \subseteq \vert B \vert \Leftrightarrow \vdash A \rightarrow B$; \\

(vii) $\vert A \vert = \vert B \vert \Leftrightarrow \vdash A \leftrightarrow B$.
\end{prop}

\begin{defn}
A estrutura $\mathcal{M} = <W, S, V>$ é um modelo canônico para $\mathbb{L}(\nabla)$,
quando satisfaz as seguintes condições: \\

- $W = \chi$; \\

- $\vert A \vert \in S(\Gamma) \Leftrightarrow \square A \in \Gamma$, para todo $\Gamma \in W$; \\

- $V(p_i) = \vert p_i \vert$, para toda variável proposicional $p_i$.
\end{defn}

\begin{prop}
Seja $\mathcal{M}$ um modelo canônico. Então, para toda fórmula A e todo $\Gamma \in W$, $\Gamma \models^{\mathcal{M}} A \Leftrightarrow A \in \Gamma$. \\
\begin{dem}
Faremos por indução sobre a complexidade das fórmulas de $\Gamma$, para o conjunto de operadores $\lbrace \neg, \land, \square \rbrace$. \\

(i) Para o caso em que $A \equiv p_i$, para cada $i \in \mathbb{N}$. Temos que $\Gamma \models^{\mathcal{M}} p_i$ see $\Gamma \in V(p_i)$ see $\Gamma \in \vert p_i \vert$. Por construção de $\vert p_i \vert$, $\Gamma$ é um conjunto em $\vert p_i \vert$ see $p_i \in \Gamma$. \\

(ii) Para o caso em que $A \equiv \neg B$ : $\Gamma \models^{\mathcal{M}} \neg B$ see $\Gamma \nvDash^{\mathcal{M}} B$. Pela hipótese de indução, temos que $B \notin \Gamma$, e como $\Gamma$ é maximalmente consistente, segue que $\neg B \in Gamma$ e, portanto, $A \in \Gamma$. \\

(iii) Para $A \equiv B \land C$ : $\Gamma \models^{\mathcal{M}} B \land C$ see $\Gamma \models^{\mathcal{M}} B$ e $\Gamma \models^{\mathcal{M}} C$. Assim, por hipótese de indução, $B, C \in \Gamma$. Pela lógica proposicional clássica, $B \land C \in \Gamma$ e, portanto, $A \in \Gamma$. \\

(iv) Para $A \equiv \square B$ : $\Gamma \models^{\mathcal{M}} \square B$ $\parallel B \parallel$ $\in S(\Gamma)$. Desse modo, pela hipótese de indução, para todo $\Omega \in W$, temos que $\Omega \models^{\mathcal{M}} B$ see $B \in \Omega$, isto é, $\parallel B \parallel$ = $\vert B \vert$. Assim, $\parallel B \parallel$ $\in S(\Gamma)$ see $\vert B \vert \in S(\Gamma)$. Agora, por definição de $S(\Gamma)$, $\vert B \vert \in S(\Gamma)$ see $\square B \in \Gamma$ e, assim, $A \in \Gamma$.
\end{dem}
\end{prop}

E agora, mais algumas definições.

\begin{defn}
Chamamos a estrutura $\mathcal{Mc}$ de menor modelo canônico quando o conjunto $S(\Gamma)$ contém apenas conjuntos de demonstrações.
\end{defn}

\begin{defn} 
A suplementação de $\mathcal{Mc}$ é o modelo $\mathcal{Mc}^{+} = <W, S^+, V>$, com a seguinte condição adicional para $S^+$. Para todo $\Gamma \in W$ e todo $X \subseteq W$: \\

$X \in S^+(\Gamma) \Leftrightarrow \exists Y \in S^+(\Gamma)$ tal que $Y \subseteq X$. 
\end{defn}

O conjunto $S^+(\Gamma)$ descrito acima pode ser facilmente reescrito da seguinte maneira $S^+(\Gamma)$ = $\lbrace X \subseteq W$ : $\vert A \vert \subseteq X$, para algum $\square A \in \Gamma \rbrace$. Assim, claramente, $S(\Gamma) \subseteq S^+(\Gamma)$. Agora, falta-nos mostrar que esta suplementação é um
modelo canônico para $\mathbb{L}(\nabla)$.

\begin{prop} 
$\mathcal{Mc}^{+}$ é um modelo canônico para $\mathbb{L}(\nabla)$. \\
\begin{dem}
Seja $\mathcal{Mc}$ o menor modelo canônico de $\mathbb{L}(\nabla)$. Para mostrarmos o que queremos, basta que verifiquemos a seguinte equivalência: para todo A e todo $\Gamma \in W, \vert A \vert \in S^+(\Gamma) \Leftrightarrow \square A \in \Gamma$. \\

$(\Rightarrow)$ Como hipótese, seja $\vert A \vert \in S^+(\Gamma)$. Daí, para algum $Y \in S(\Gamma)$, temos que $Y \subseteq \vert A \vert$. Como $\mathcal{Mc}$ é o menor modelo canônico, então $Y = \vert B \vert$, para algum B. Desse modo, $\vert B \vert \subseteq \vert A \vert$ e $\square B \in \Gamma$. Assim, pelos resultados anteriores, segue que $\vdash A \rightarrow B$ e, agora, por RE temos que, $\vdash \square A \rightarrow \square B$. Portanto, $\square A \in \Gamma$. \\

$(\Leftarrow)$ Como hipótese, agora, consideremos $\square A \in \Gamma$. Daí, $\vert A \vert \in S(\Gamma)$, e como $S(\Gamma) \subseteq S^+(\Gamma)$, vem que $\vert A \vert \in S^+(\Gamma)$.
\end{dem}
\end{prop}

\begin{prop}
Seja $\mathcal{Mc}$ o menor modelo canônico para $\mathbb{L}(\nabla)$ e $\mathcal{Mc}^{+}$ sua suplementação. Então as condições (c), (h), (t) e (n) são válidas em $\mathcal{Mc}^{+}$. \\
\begin{dem}
Pelo resultado anterior, segue que $\mathcal{Mc}^{+}$ é um modelo canônico para $\mathbb{L}(\nabla)$. Assim, basta verificarmos as condições listadas.\\

- Para (c): Sejam $\Gamma \in W$ e X e Y subconjuntos de W tais que $X \in S^+(\Gamma)$ e $Y \in S^+(\Gamma)$. Pela definição de $S^+(\Gamma)$, segue que existe $Z \subseteq X$ tal que $Z \in S(\Gamma)$ e também que existe $U \subseteq Y$ tal que $U \in S(\Gamma)$. Assim, pela condição (c), segue que $Z \cap U \in S(\Gamma)$ e, novamente pela definição de $S^+(\Gamma)$, vem que $X \cap Y \in S^+(\Gamma)$. Portanto,
(c) se mantém em $\mathcal{Mc}^{+}$. \\

- Para (h): Sejam $\Gamma \in W$ e X e Y subconjuntos de W tais que $X \in S^+(\Gamma)$ ou $Y \in S^+(\Gamma)$. Por definição de $S^+(\Gamma)$, segue que existe $Z \subseteq X$ tal que $Z \in S(\Gamma)$ ou que existe $U \subseteq Y$ tal que $U \in S(\Gamma$). Assim, pela condição (h), segue que $Z \cup U \in S(\Gamma)$ e, pela definição de $S^+(\Gamma)$, $X \cup Y \in S^+(\Gamma)$. Portanto, (h) se mantém em $\mathcal{Mc}^{+}$. \\

- Para (t): Sejam $\Gamma \in W$ e X subconjunto de W tal que $X \in S^+(\Gamma)$. Pela definição de $S^+(\Gamma)$, existe $Y \subseteq X$ tal que $Y \in S(\Gamma)$. Assim, por (t), $\Gamma \in Y$ e como $Y \subseteq X$, então $\Gamma \in X$. Portanto, (t) se mantém em $\mathcal{Mc}^{+}$. \\

- Para (n): Queremos mostrar que $W \in S^+(\Gamma)$. Para isso, consideremos um conjunto de fórmulas $\Gamma \in W$ e que, de fato, $W \subseteq W$. Desde que $W = \chi$, e $\Gamma \in W$, para algum $A \in \Gamma$, temos que $\square A \in \Gamma$, e daí, $\vert A \vert \in W$. Desse modo, pela definição de $S^+(\Gamma)$, segue que $W \in S^+(\Gamma)$. Portanto, (n) se mantém em $\mathcal{Mc}^{+}$.
\end{dem}
\end{prop}

Agora, mais uma definição e mostraremos o que queremos.

\begin{defn}
Seja $\mathcal{Mc}$ o menor modelo canônico e $\mathcal{Mc}^{+}$ sua suplementação. O modelo $\mathcal{Mc}^{+}_{\Gamma} = <W_{\Gamma}, S_{\Gamma}, V_{\Gamma}>$ será denominado suplementado para $\Gamma$ quando as
seguintes condições são válidas na estrutura: \\

- $W_{\Gamma}$ = $\lbrace \Sigma \in \chi$ : $\Gamma \subseteq \Sigma \rbrace$; \\

- Para todo $\Sigma \in W_{\Gamma}$ e todo $X \subseteq W_{\Gamma}, X \in S_{\Gamma} \Leftrightarrow$ para algum $Y \in S^+(\Sigma)$, $X = Y \cap W_{\Gamma}$; \\

- $V_{\Gamma}(p_i) = \lbrace \Sigma \in \chi$ : $p_i \in \Sigma$ e $ \Gamma \in \Sigma \rbrace$, para toda variável proposicional $p_i$.
\end{defn}

Intuitivamente, a estrutura acima é obtida excluindo de $\mathcal{Mc}^{+}$ todos os conjuntos maximalmente consistentes que não contém $\Gamma$. Dessa forma, $\mathcal{Mc}^{+}_{\Gamma}$ é um modelo canônico e, portanto, as condições (c), (h), (t) e (n) continuam válidas nesta estrutura. E mais, construído dessa forma, para todo $\Sigma \in W_{\Gamma}, \Gamma \subseteq \Sigma$, segue que $\Sigma \models^{\mathcal{K}} \Gamma$, com $\mathcal{K} = \mathcal{Mc}^{+}_{\Gamma}$ e, portanto, $\models^{\mathcal{K}} \Gamma$.

Agora, finalmente, mostraremos com o auxílio destas estruturas a completude.

\begin{prop} $\textit{(Completude)}$
Se $\Gamma \models A$, então $\Gamma \vdash A$. \\
\begin{dem}
Façamos pela recíproca da contrária. Assim, assumamos que $\Gamma \nvdash A$. Desse modo, obtemos que $\Gamma \nvdash \neg \neg A$ e, então, que $\Gamma \cup \lbrace \neg A \rbrace$ é consistente. Daí, pelo Teorema de Lindenbaum, existe um conjunto maximalmente consistente, $\Delta$, que estende $\Gamma \cup \lbrace \neg A \rbrace$, i.e., $\Gamma \cup \lbrace \neg A \rbrace \subseteq \Delta$, com $\neg A \in \Delta$ e $A \notin \Delta$. Agora, $\Delta$ é um conjunto maximalmente consistente em $\mathbb{L}(\nabla)$ e $\Gamma \subseteq \Delta$. Logo, considerando a estrutura $\mathcal{Mc}^{+}_{\Gamma}$, temos que $\Delta$ é um conjunto tal que para o modelo $\mathcal{Mc}^{+}_{\Gamma}$ a fórmula A falha. Dessa maneira, $\nvDash^{\mathcal{K}} A$, para $\mathcal{K} = \mathcal{Mc}^{+}_{\Gamma}$ e, portanto, $\nvDash A$.
\end{dem}
\end{prop}

Portanto, o presente trabalho apresenta uma semântica de vizinhança que é equivalente à versão axiomática da lógica proposicional do plausível.


\section*{Considerações finais}

${}$ \hspace {0,5 cm} O presente trabalho mostra uma semântica relacional para a lógica proposicional do plausível. Inicialmente, precisamos considerar um novo sistema modal, o qual denotamos por $\mathbb{LP}(\square)$, e denominamos de lógica proposicional modal do plausível, de maneira que as ideias intuitivas do operador do plausível, $\nabla$, pudessem ser capturadas pelo operador de necessidade, $\square$. Isto é, mostramos a equivalência dedutiva entre essas duas abordagens, para podermos tratar com conceitos da literatura sobre este último operador.

Feito isso, passamos a pensar sobre as semânticas de Kripke. No entanto, as semânticas de Kripke usuais são aplicadas em lógicas modais normais, ou seja, aquelas em que o axioma $\textbf{K}$ se verifica. Dado que a lógica proposicional do plausível é uma lógica subnormal, precisamos de um outro ambiente para obtermos uma semântica de caráter relacional.

E assim optamos pelas semânticas relacionais ou de vizinhança. Pois as mesmas tratam de certa forma dos mundos possíveis e, no nosso caso, de algumas relações especiais impostas nas definições dos modelos que investigamos.

Assim, para determinarmos quais seriam as relações impostas, acabamos por olhar para os axiomas e obtivemos algumas relações que de algum modo os espelhassem. Ao apresentarmos tais relações, mostramos a correção do sistema lógico com esta nova semântica. 

Já para mostrarmos a completude tivemos que seguir um caminho que é usual na literatura. Passamos a considerar os conjuntos maximalmente consistentes. E mais, usamos de um resultado, do ponto de vista teórico, muito forte que é o teorema de Lindenbaum. Consideramos ser este um resultado extremamente interessante, pois ele é, de algum modo, equivalente ao axioma da escolha.

Portanto, podemos dizer que o presente artigo pretendia mostrar uma nova semântica para a lógica proposicional do plausível e este objetivo foi alcançado.


\section*{Agradecimentos}

Agradecemos apoio da FAPESP e do DM da UNESP - Câmpus de Bauru.


\end{document}